\documentclass[brochure,12pt]{bourbaki}
\usepackage[all]{xy}
\usepackage{amssymb,amsfonts,amsmath}
\usepackage[T1]{fontenc}
\newcommand{\og}{\guillemotleft\hspace*{1pt}}
\newcommand{\fg}{\hspace*{1pt}\guillemotright}

\newcommand{\ie}{{\em i.e.}\ }

\newcommand{\ko}{\: , \;}

\newcommand{\ol}[1]{\overline{#1}}

\numberwithin{equation}{subsection}
\newtheorem{theorem}{Theorem}
\numberwithin{theorem}{section}
\newtheorem{classification-theorem}[theorem]{Classification Theorem}
\newtheorem{decomposition-theorem}[theorem]{Decomposition Theorem}
\newtheorem{proposition-definition}[theorem]{Proposition-Definition}
\newtheorem{periodicity-conjecture}[theorem]{Periodicity Conjecture}

\newtheorem{proposition}[theorem]{Proposition}

\newcommand{\reminder}[1]{}

\newcommand{\opname}[1]{\operatorname{\mathsf{#1}}}

\newcommand{\rep}{\opname{rep}\nolimits}

\renewcommand{\rep}{\opname{rep}\nolimits}

\newcommand{\scr}{\mathcal}

\newcommand{\Z}{\mathbb{Z}}
\newcommand{\N}{\mathbb{N}}
\newcommand{\Q}{\mathbb{Q}}
\newcommand{\C}{\mathbb{C}}

\newcommand{\la}{\leftarrow}
\newcommand{\iso}{\stackrel{_\sim}{\rightarrow}}

\newcommand{\Hom}{\opname{Hom}}

\newcommand{\Ext}{\opname{Ext}}

\newcommand{\ten}{\otimes}

\newcommand{\ca}{{\mathcal A}}

\newcommand{\cc}{{\mathcal C}}
\newcommand{\cd}{{\mathcal D}}

\newcommand{\cF}{{\mathcal F}}
\newcommand{\cg}{{\mathcal G}}

\newcommand{\cm}{{\mathcal M}}
\newcommand{\cn}{{\mathcal N}}

\newcommand{\cp}{{\mathcal P}}

\renewcommand{\phi}{\varphi}

\renewcommand{\hat}[1]{\widehat{#1}}

\newcommand{\sgn}{\mbox{sgn}}

\newcommand{\HH}{HH}

\renewcommand{\tilde}[1]{\widetilde{#1}}

\newcommand{\arr}[1]{\stackrel{#1}{\rightarrow}}
\newcommand{\maparr}[1]{\stackrel{#1}{\mapsto}}

\date{Novembre 2009}
\bbkannee{62\`eme ann\'ee, 2009-2010}
\bbknumero{1014}
\title{ALG\`EBRES AMASS\'EES ET APPLICATIONS}
\subtitle{d'apr\`es Fomin-Zelevinsky, \ldots}
\author{Bernhard KELLER}
\address{Universit\'e Paris Diderot -- Paris 7\\
Institut de Math\'ematiques de Jussieu\\
U.M.R.~7586 du CNRS\\
U.F.R.~de Math\'ematiques\\
 Case 7012\\
B\^atiment Chevaleret\\
F--75205 Paris Cedex 13}
\email{keller@math.jussieu.fr}

\begin{document}
\maketitle

% \tableofcontents

% \og Hello\fg

\section*{INTRODUCTION}
% Sergey Fomin et Andrei Zelevinsky ont invent\'e
% \cite{FominZelevinsky02} les alg\`ebres amass\'ees (cluster algebras)
% au d\'ebut des ann\'ees 2000 dans le but de fournir un cadre
% combinatoire \`a l'\'etude des bases canoniques dans les groupes
% quantiques \cite{Kashiwara90} \cite{Lusztig90} et de la positivit\'e
% totale dans les groupes alg\'ebriques \cite{Lusztig96}. En partie en
% collaboration avec Arkady Berenstein, ils les ont d\'evelopp\'ees dans
%  une s\'erie d'articles \cite{FominZelevinsky03}
% \cite{BerensteinFominZelevinsky05} \cite{FominZelevinsky07}.
Les alg\`ebres amass\'ees (cluster algebras), invent\'ees \cite{FominZelevinsky02} par
Sergey Fomin et Andrei Zelevinsky au d\'ebut des ann\'ees 2000,
sont des alg\`ebres commutatives, dont les g\'en\'erateurs
et les relations sont construits de fa\c{c}on r\'ecursive.
Parmi ces alg\`ebres se trouvent les alg\`ebres de coordonn\'ees homog\`enes
sur les grassmanniennes, les vari\'et\'es de drapeaux et
beaucoup d'autres vari\'et\'es qui jouent un r\^ole important en g\'eom\'etrie
et th\'eorie des repr\'esentations. La motivation principale de
Fomin et Zelevinsky \'etait de trouver un cadre combinatoire
pour l'\'etude des bases canoniques dont on dispose
\cite{Kashiwara90} \cite{Lusztig90} dans ces alg\`ebres
et qui sont \'etroitement li\'ees \`a
la notion de positivit\'e totale \cite{Lusztig96}
dans les vari\'et\'es associ\'ees.
Il s'est
av\'er\'e rapidement que la combinatoire des alg\`ebres amass\'ees
intervenait \'egalement dans de nombreux autres sujets, par exemple
dans
\begin{itemize}

\item la g\'eom\'etrie de Poisson
\cite{GekhtmanShapiroVainshtein03}
\cite{GekhtmanShapiroVainshtein05}
\cite{GekhtmanShapiroVainshtein08}
\cite{BerensteinZelevinsky05} \ldots ;

\item les syst\`emes dynamiques discrets
\cite{FominZelevinsky03b}
\cite{Kedem08}
\cite{DiFrancescoKedem08}
\cite{InoueIyamaKunibaNakanishiSuzuki08} \ldots ;

\item les espaces de Teichm\"uller sup\'erieurs
\cite{FockGoncharov03}
\cite{FockGoncharov07b}
\cite{FockGoncharov06}
\cite{FockGoncharov09}
\ldots ;

\item la combinatoire et en particulier l'\'etude de
poly\`edres tels les associa\`edres de Sta\-sheff
\cite{ChapotonFominZelevinsky02}
\cite{Chapoton04}
\cite{IngallsThomas06}
\cite{Krattenthaler06}
\cite{FominReading05}
\cite{FominShapiroThurston08}
\cite{Musiker07}
\cite{MusikerSchifflerWilliams09}
\ldots ;

\item la g\'eom\'etrie alg\'ebrique (commutative ou
non commutative) et en particulier l'\'etude des
conditions de stabilit\'e de Bridgeland
% \cite{Bridgeland06a}
% \cite{Bridgeland06}
\cite{Bridgeland07},
les alg\`ebres Calabi-Yau
\cite{IyamaReiten08} \cite{Ginzburg06},
les invariants de Donaldson-Thomas
% \cite{Szendroi07}
% \cite{Kontsevich07a} \cite{Kontsevich07}
\cite{JoyceSong09}
\cite{KontsevichSoibelman08}
\cite{Reineke09}
\cite{GaiottoMooreNeitzke09}
\ldots ;

\item et la th\'eorie des repr\'esentations des carquois
et des alg\`ebres de dimension finie, voir par exemple
les articles de synth\`ese
\cite{BuanMarsh06} \cite{Reiten06} \cite{Ringel07}
\cite{GeissLeclercSchroeer08a}
\cite{Keller08c}.
\end{itemize}
Nous renvoyons aux articles d'initiation
\cite{FominZelevinsky03a}
\cite{Zelevinsky04}
\cite{Zelevinsky02}
\cite{Zelevinsky05}
\cite{Zelevinsky07a}
et au portail des alg\`ebres amass\'ees \cite{Fomin07}
pour plus d'informations sur les alg\`ebres amass\'ees et
leurs liens avec d'autres sujets math\'ematiques
(et physiques).

Dans cet expos\'e, nous donnons une introduction concise
aux alg\`ebres amass\'ees (section~1) et pr\'esentons deux applications~:
\begin{itemize}

\item la d\'emonstration de la p\'eriodicit\'e de certains
syst\`emes dynamiques discrets, d'apr\`es
Fomin-Zelevinsky \cite{FominZelevinsky03b} et l'auteur \cite{Keller08c} \cite{Keller10a}
(section~\ref{ss:Y-graines})~;

\item la construction de bases duales semi-canoniques,
d'apr\`es Geiss-Leclerc-Schr\"oer \cite{GeissLeclercSchroeer06}
(section~\ref{ss:application-aux-bases-duales-semi-canoniques}).

\end{itemize}
Ces applications sont fond\'ees sur la cat\'egorification additive des
alg\`ebres amass\'ees \`a l'aide de cat\'egories de repr\'esentations
de carquois (avec relations).  Nous en d\'ecrivons les id\'ees
principales \`a la section~\ref{s:categorifications}. Nous y
esquissons \'egalement des d\'eveloppements r\'ecents importants li\'es \`a la
cat\'egorification mono\"{\i}dale d'alg\`ebres amass\'ees
\cite{HernandezLeclerc09} \cite{Nakajima09} et \`a leur \'etude via
les carquois \`a potentiel \cite{DerksenWeymanZelevinsky08}
\cite{DerksenWeymanZelevinsky09}.

\section{DESCRIPTION ET PREMIERS EXEMPLES}
\subsection{Description}
Une {\em alg\`ebre amass\'ee} est une $\Q$-alg\`ebre commutative munie
d'un ensemble de g\'en\'erateurs distingu\'es (les {\em variables d'amas})
regroup\'es dans des parties (les {\em amas}) de cardinal constant
(le {\em rang}) qui sont construites r\'ecursivement par {\em mutation}
\`a partir d'un {\em amas initial}.
L'ensemble des variables d'amas peut \^etre fini ou infini.

\begin{theo}[\cite{FominZelevinsky03}] Les alg\`ebres amass\'ees
n'ayant qu'un nombre fini de variables d'amas sont param\'etr\'ees par
les syst\`emes de racines finis.
\end{theo}

La classification est donc analogue \`a celle des alg\`ebres
de Lie semi-simples complexes.
Nous allons pr\'eciser le th\'eor\`eme (dans le cas simplement lac\'e)
\`a la section~\ref{s:algebres-amassees-associees-aux-carquois}.

\subsection{Premier exemple}
Pour illustrer la description et le th\'eor\`eme, pr\'esentons \cite{Zelevinsky07a}
l'alg\`ebre amass\'ee $\ca_{A_2}$ associ\'ee au syst\`eme de racines $A_2$.
Par d\'efinition, elle est engendr\'ee sur $\Q$ par les variables
d'amas $x_m$, $m\in \Z$, soumises aux {\em relations d'\'echange}
\[
x_{m-1} x_{m+1} = 1+x_m \ko m\in \Z.
\]
Ses amas sont par d\'efinition les paires de variables cons\'ecutives
$\{ x_m, x_{m+1}\}$, $m\in \Z$. L'amas initial est $\{x_1, x_2\}$ et
deux amas sont reli\'es par une mutation si et seulement si ils ont
exactement une variable d'amas en commun.

Les relations d'\'echange permettent d'exprimer toute variable d'amas
comme fonction rationnelle des variables initiales $x_1, x_2$ et donc
d'identifier l'alg\`ebre $\ca_{A_2}$ \`a une sous-alg\`ebre du corps
$\Q(x_1, x_2)$. Afin d'expliciter cette sous-alg\`ebre, calculons les
$x_m$ pour $m\geq 3$.  Nous avons~:
\begin{align}
x_3 &= \frac{1+x_2}{x_1} \\
x_4 &= \frac{1+x_3}{x_2} = \frac{x_1 + 1 + x_2}{x_1 x_2}\\
x_5 &= \frac{1+x_4}{x_3}
= \frac{x_1 x_2 + x_1 +1 +x_2}{x_1 x_2}\; \div \; \frac{1+x_2}{x_1}
= \frac{1+x_1}{x_2} \; . \label{eq:laurent}
\end{align}
Notons que, contrairement \`a ce qu'on pourrait attendre,
le d\'enominateur dans~\ref{eq:laurent} reste un mon\^ome !
En fait, toute variable d'amas dans une alg\`ebre amass\'ee quelconque
est un polyn\^ome de Laurent, voir le th\'eor\`eme~\ref{thm:class-finie-detaillee}.
 Continuons le calcul~:
\begin{align}
x_6 &= \frac{1+x_5}{x_4}
     =  \frac{x_2+1+x_1}{x_2}\; \div \;\frac{x_1+1+x_2}{x_1 x_2} = x_1 \\
x_7 &= (1+x_1)\; \div \; \frac{1+x_1}{x_2}  = x_2.
\end{align}
Il est alors clair que la suite des $x_m$, $m\in\Z$, est $5$-p\'eriodique
et que le nombre de variables d'amas est effectivement fini
et \'egal \`a cinq. Outre les deux variables initiales
$x_1$ et $x_2$ nous avons trois variables non initiales
$x_3$, $x_4$ et $x_5$. En examinant leurs d\'enominateurs,
nous voyons qu'elles sont en bijection naturelle
avec les racines positives $\alpha_1$, $\alpha_1+\alpha_2$, $\alpha_2$
du syst\`eme de racines de type $A_2$. Ceci se g\'en\'eralise
\`a tout diagramme de Dynkin, voir le th\'eor\`eme~\ref{thm:class-finie-detaillee}.

\subsection{Alg\`ebres amass\'ees de rang $2$} \`A tout couple d'entiers
positifs $(b,c)$ est associ\'ee  une alg\`ebre amass\'ee $\ca_{(b,c)}$. On la
d\'efinit de la m\^eme mani\`ere que $\ca_{A_2}$, mais en rempla\c{c}ant les relations d'\'echange par
\[
x_{m-1} x_{m+1} = \left\{ \begin{array}{ll} x_m^b +1 & \mbox{si $m$ est impair, } \\
                                            x_m^c +1 & \mbox{si $m$ est pair.} \end{array} \right.
\]
L'alg\`ebre $\ca_{(b,c)}$ n'a qu'un nombre fini de variables d'amas
si et seulement si $bc\leq 3$, autrement dit si la matrice
\[
\left[ \begin{array}{cc} 2 & -b \\ -c & 2 \end{array} \right]
\]
est la matrice de Cartan d'un syst\`eme de racines $\Phi$ de rang $2$.
Le lecteur pourra s'amuser \`a v\'erifier que, dans ce cas, les
variables d'amas non initiales sont toujours param\'etr\'ees par
les racines positives de $\Phi$.

\section{ALG\`EBRES AMASS\'EES ASSOCI\'EES AUX CARQUOIS}
\label{s:algebres-amassees-associees-aux-carquois}

\subsection{Mutation des carquois} Un {\em carquois} est un graphe orient\'e, c'est-\`a-dire
un quadruplet $Q=(Q_0, Q_1, s,t)$ form\'e d'un ensemble de sommets $Q_0$, d'un
ensemble de fl\`eches $Q_1$ et de deux applications $s$ et $t$ de $Q_1$ dans $Q_0$
qui, \`a une fl\`eche $\alpha$, associent respectivement sa \underline{s}ource et son
bu\underline{t}. En pratique, on repr\'esente un carquois par un dessin comme dans
l'exemple qui suit~:
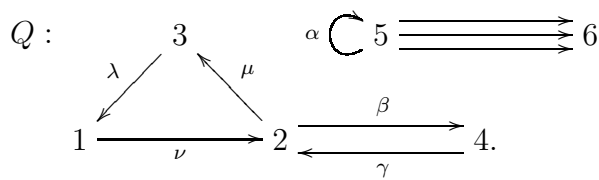
\[ Q:
\xymatrix{ & 3 \ar[ld]_\lambda & & 5 \ar@(dl,ul)[]^\alpha \ar@<1ex>[rr] \ar[rr] \ar@<-1ex>[rr] & & 6 \\
  1 \ar[rr]_\nu & & 2 \ar@<1ex>[rr]^\beta \ar[ul]_\mu & & 4.
  \ar@<1ex>[ll]^\gamma }
\]
Une fl\`eche $\alpha$ dont la source
et le but co\"{\i}ncident est une {\em boucle}~;
un {\em $2$-cycle} est un couple de fl\`eches distinctes
$\beta$ et $\gamma$ telles que $s(\beta)=t(\gamma)$ et $t(\beta)=s(\gamma)$.
De m\^eme, on d\'efinit les $n$-cycles pour tout entier positif $n$.
Un sommet $i$ d'un carquois est une {\em source} (respectivement un {\em puits})
s'il n'existe aucune fl\`eche de but $i$ (respectivement de source~$i$).

Appelons {\em bon carquois} un carquois fini sans boucles
ni $2$-cycles dont l'ensemble des sommets est l'ensemble des entiers
$1$, \ldots, $n$ pour un entier positif $n$. \`A un isomorphisme
fixant les sommets pr\`es, un tel carquois $Q$ est donn\'e par la
matrice antisym\'etrique $B=B_Q$ dont le coefficient $b_{ij}$ est
la diff\'erence entre le nombre de fl\`eches de $i$ \`a $j$ et
le nombre de fl\`eches de $j$ \`a $i$ pour tous $1\leq i,j\leq n$.
R\'eciproquement, toute matrice antisym\'etrique $B$ \`a coefficients
entiers provient d'un bon carquois $Q_B$.

Soient $Q$ un bon carquois et $k$ un sommet de $Q$. La {\em carquois mut\'e $\mu_k(Q)$}
est le carquois obtenu \`a partir de $Q$ comme suit~:
\begin{itemize}
\item[(1)] pour tout sous-carquois $\xymatrix{i \ar[r]^\beta & k \ar[r]^\alpha & j}$,
on rajoute une nouvelle fl\`eche $[\alpha\beta]: i \to j$~;
\item[(2)] on renverse toutes les fl\`eches de source ou de but $k$~;
\item[(3)] on supprime les fl\`eches d'un ensemble maximal de $2$-cycles disjoints
deux \`a deux.
\end{itemize}
Si $B$ est la matrice antisym\'etrique associ\'ee \`a $Q$ et $B'$ celle associ\'ee
\`a $\mu_k(Q)$, on a
\[
b'_{ij} = \left\{ \begin{array}{ll}
-b_{ij} & \mbox{si $i=k$ ou $j=k$~;} \\
b_{ij}+\sgn(b_{ik})\max(0, b_{ik} b_{kj}) & \mbox{sinon.}
\end{array} \right.
\]
C'est la r\`egle de {\em mutation des matrices} antisym\'etriques
(plus g\'en\'eralement~: antisym\'etrisables) introduite
par Fomin-Zelevinsky dans \cite{FominZelevinsky02},
voir aussi \cite{FominZelevinsky07}.

On v\'erifie sans peine que $\mu_k$ est une involution. Par
exemple, les carquois
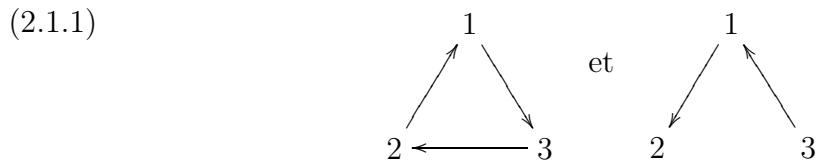
\begin{equation} \label{quiver1}
\begin{xy} 0;<0.3pt,0pt>:<0pt,-0.3pt>::
(94,0) *+{1} ="0",
(0,156) *+{2} ="1",
(188,156) *+{3} ="2",
"1", {\ar"0"},
"0", {\ar"2"},
"2", {\ar"1"},
\end{xy}
\begin{minipage}{1cm}
\vspace*{1cm}
\begin{center} et \end{center}
\end{minipage}
\begin{xy} 0;<0.3pt,0pt>:<0pt,-0.3pt>::
(92,0) *+{1} ="0",
(0,155) *+{2} ="1",
(188,155) *+{3} ="2",
"0", {\ar"1"},
"2", {\ar"0"},
\end{xy}
\end{equation}
sont reli\'es par la mutation par rapport au sommet $1$. Notons
que, du point de la th\'eorie des
repr\'esentations, ces carquois sont tr\`es diff\'erents.
Deux carquois sont {\em \'equivalents par mutation} s'ils sont reli\'es
par une suite finie de mutations.
On v\'erifie facilement, par exemple \`a l'aide de
\cite{KellerQuiverMutationApplet}, que les trois carquois suivants sont
\'equivalents par mutation
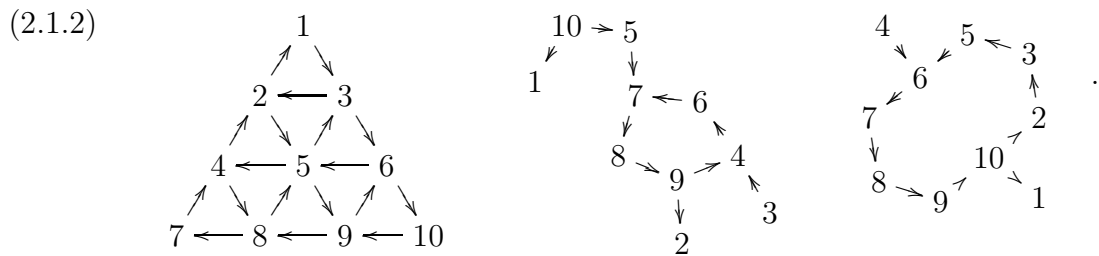
\begin{equation} \label{quiver3}
\begin{xy} 0;<0.6pt,0pt>:<0pt,-0.6pt>::
(79,0) *+{1} ="0",
(52,44) *+{2} ="1",
(105,44) *+{3} ="2",
(26,88) *+{4} ="3",
(79,88) *+{5} ="4",
(131,88) *+{6} ="5",
(0,132) *+{7} ="6",
(52,132) *+{8} ="7",
(105,132) *+{9} ="8",
(157,132) *+{10} ="9",
"1", {\ar"0"},
"0", {\ar"2"},
"2", {\ar"1"},
"3", {\ar"1"},
"1", {\ar"4"},
"4", {\ar"2"},
"2", {\ar"5"},
"4", {\ar"3"},
"6", {\ar"3"},
"3", {\ar"7"},
"5", {\ar"4"},
"7", {\ar"4"},
"4", {\ar"8"},
"8", {\ar"5"},
"5", {\ar"9"},
"7", {\ar"6"},
"8", {\ar"7"},
"9", {\ar"8"},
\end{xy}
\quad\quad
\begin{xy} 0;<0.3pt,0pt>:<0pt,-0.3pt>::
(0,70) *+{1} ="0",
(183,274) *+{2} ="1",
(293,235) *+{3} ="2",
(253,164) *+{4} ="3",
(119,8) *+{5} ="4",
(206,96) *+{6} ="5",
(125,88) *+{7} ="6",
(104,164) *+{8} ="7",
(177,194) *+{9} ="8",
(39,0) *+{10} ="9",
"9", {\ar"0"},
"8", {\ar"1"},
"2", {\ar"3"},
"3", {\ar"5"},
"8", {\ar"3"},
"4", {\ar"6"},
"9", {\ar"4"},
"5", {\ar"6"},
"6", {\ar"7"},
"7", {\ar"8"},
\end{xy}
\quad\quad
\begin{xy} 0;<0.3pt,0pt>:<0pt,-0.3pt>::
(212,217) *+{1} ="0",
(212,116) *+{2} ="1",
(200,36) *+{3} ="2",
(17,0) *+{4} ="3",
(123,11) *+{5} ="4",
(64,66) *+{6} ="5",
(0,116) *+{7} ="6",
(12,196) *+{8} ="7",
(89,221) *+{9} ="8",
(149,166) *+{10} ="9",
"9", {\ar"0"},
"1", {\ar"2"},
"9", {\ar"1"},
"2", {\ar"4"},
"3", {\ar"5"},
"4", {\ar"5"},
"5", {\ar"6"},
"6", {\ar"7"},
"7", {\ar"8"},
"8", {\ar"9"},
\end{xy}
\begin{minipage}{1cm}
\vspace*{1.5cm}
\begin{center} . \end{center}
\end{minipage}
\end{equation}
La classe de mutation commune de ces carquois comporte 5739 carquois (\`a
isomorphisme pr\`es). La classe de mutation de la \og plupart\fg\ des
carquois est infinie.  La classification des carquois ayant une
classe de mutation finie est un probl\`eme difficile, r\'esolu
r\'ecemment par Felikson-Shapiro-Tumarkin
\cite{FeliksonShapiroTumarkin08}~: outre les carquois \`a deux
sommets et les carquois associ\'es \`a des surfaces \`a bord
marqu\'ees \cite{FominShapiroThurston08}, la liste contient $11$
carquois exceptionnels, dont le plus grand est dans la classe de
mutation des carquois~\ref{quiver3}.

\subsection{Mutation des graines, alg\`ebres amass\'ees}
Soient $n\geq 1$ un entier et $\cF$ le corps $\Q(x_1, \ldots, x_n)$
engendr\'e par $n$ ind\'etermin\'ees $x_1, \ldots, x_n$.
Une {\em graine} (appel\'ee aussi {\em $X$-graine}) est un couple $(R,u)$,
o\`u $R$ est un bon carquois et $u$ une suite $u_1, \ldots, u_n$ qui engendre
librement le corps $\cF$. Si $(R,u)$ est une graine et
$k$ un sommet de $R$, la {\em mutation $\mu_k(R,u)$} est
la graine $(R',u')$, o\`u $R'=\mu_k(R)$ et $u'$ est
obtenu \`a partir de $u$ en rempla\c{c}ant l'\'el\'ement
$u_k$ par l'\'el\'ement $u_k'$ d\'efini par la
{\em relation d'\'echange}
\begin{equation} \label{eq:echange}
u_k' u_k= \prod_{s(\alpha)=k} u_{t(\alpha)} + \prod_{t(\alpha)=k} u_{s(\alpha)} .
\end{equation}
On v\'erifie que $\mu_k^2(R,u)=(R,u)$. Par exemple, les mutations de la graine
\[
\xymatrix{1 \ar[r] & 2 \ar[r] & 3} \ko \{x_1, x_2, x_3\}
\]
par rapport aux sommets $1$ et $2$ sont les graines
\begin{equation} \label{eq:seeds-for-A3}
\xymatrix{1 & 2 \ar[l] \ar[r] & 3} \ko \{\frac{1+x_2}{x_1}, x_2, x_3\}
\quad\mbox{ et }\quad
\xymatrix{1 \ar@/^1pc/[rr] & 2 \ar[l] & 3 \ar[l]} \ko \{ x_1, \frac{x_1+x_3}{x_2}, x_3\}.
\end{equation}

Fixons maintenant un bon carquois $Q$. La {\em graine initiale} est
$(Q,\{x_1, \ldots, x_n\})$. Un {\em amas} associ\'e \`a $Q$ est une
suite $u$ qui appara\^{\i}t dans une graine $(R,u)$ obtenue \`a partir
de la graine initiale par mutation it\'er\'ee. Les {\em variables d'amas}
sont les \'el\'ements des amas. {\em L'alg\`ebre amass\'ee $\ca_Q$}
est la sous-alg\`ebre de $\cF$ engendr\'ee par les variables d'amas.
Clairement, si $(R,u)$ est une graine associ\'ee \`a $Q$, l'isomorphisme
naturel
\[
\Q(u_1, \ldots, u_n) \iso \Q(x_1, \ldots, x_n)
\]
induit un isomorphisme de $\ca_R$ sur $\ca_Q$ qui pr\'eserve
les variables d'amas et les amas. L'alg\`ebre amass\'ee $\ca_Q$ est
donc un invariant de la classe de mutation de $Q$. Il est utile
d'introduire un objet combinatoire qui code la construction
r\'ecursive des graines~: le {\em graphe d'\'echange}. Par d\'efinition,
ses sommets sont les classes d'isomorphisme de graines (les isomorphismes
renum\'erotent les sommets et les variables) et ses ar\^etes
correspondent aux mutations. Par exemple, le graphe d'\'echange
obtenu \`a partir du carquois
$
Q: \xymatrix{1 \ar[r] & 2 \ar[r] & 3}
$
est le $1$-squelette de l'associa\`edre de Stasheff
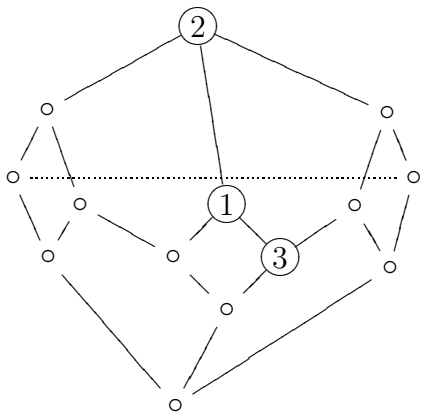
\cite{Stasheff63}
\cite{ChapotonFominZelevinsky02}:
\[
\begin{xy} 0;<0.4pt,0pt>:<0pt,-0.4pt>::
(173,0) *+<8pt>[o][F]{2} ="0",
(0,143) *+{\circ} ="1",
(63,168) *+{\circ} ="2",
(150,218) *+{\circ} ="3",
(250,218) *+<8pt>[o][F]{3} ="4",
(375,143) *+{\circ} ="5",
(350,82) *+{\circ} ="6",
(152,358) *+{\circ} ="7",
(200,168) *+<8pt>[o][F]{1} ="8",
(200,268) *+{\circ} ="9",
(32,79) *+{\circ} ="10",
(33,218) *+{\circ} ="11",
(320,170) *+{\circ} ="12",
(353,228) *+{\circ} ="13",
"0", {\ar@{-}"6"},
"0", {\ar@{-}"8"},
"0", {\ar@{-}"10"},
"1", {\ar@{.}"5"},
"1", {\ar@{-}"10"},
"11", {\ar@{-}"1"},
"2", {\ar@{-}"3"},
"10", {\ar@{-}"2"},
"2", {\ar@{-}"11"},
"3", {\ar@{-}"8"},
"9", {\ar@{-}"3"},
"8", {\ar@{-}"4"},
"4", {\ar@{-}"9"},
"4", {\ar@{-}"12"},
"6", {\ar@{-}"5"},
"5", {\ar@{-}"13"},
"12", {\ar@{-}"6"},
"9", {\ar@{-}"7"},
"11", {\ar@{-}"7"},
"13", {\ar@{-}"7"},
"13", {\ar@{-}"12"},
\end{xy}
\]
Le sommet $1$ correspond \`a la graine initiale et les sommets
$2$ et $3$ aux graines \ref{eq:seeds-for-A3}.

Fixons un bon carquois connexe $Q$. Si son graphe sous-jacent est
un diagramme de Dynkin simplement lac\'e de type $\Delta$, nous
disons que $Q$ est un {\em carquois de Dynkin de type $\Delta$.}

\begin{theo}[\cite{FominZelevinsky03}] \label{thm:class-finie-detaillee}
\begin{itemize}
\item[(a)] Toute variable d'amas de $\ca_Q$ est un polyn\^ome de Laurent \`a coefficients entiers
\cite{FominZelevinsky02}.
\item[(b)] L'alg\`ebre amass\'ee $\ca_Q$ n'a qu'un nombre fini de variables d'amas si et seulement
si $Q$ est \'equivalent par mutation \`a un carquois de Dynkin.
Dans ce cas, le graphe~$\Delta$ sous-jacent \`a ce carquois est unique et s'appelle
le {\em type amass\'e de $Q$}.
\item[(c)] Si $Q$ est un carquois de Dynkin de type $\Delta$, alors les variables
d'amas non initiales de $\ca_Q$ sont en bijection avec les racines positives du syst\`eme de racines
$\Phi$ de $\Delta$~; plus pr\'ecis\'ement, si $\alpha_1$, \ldots, $\alpha_n$ sont
les racines simples, alors pour toute racine positive $\alpha=d_1 \alpha_1 + \cdots + d_n \alpha_n$,
il existe une unique variable d'amas non initiale $X_\alpha$ de d\'enominateur
$x_1^{d_1} \cdots  x_n^{d_n}$.
\end{itemize}
\end{theo}

Un {\em mon\^ome d'amas} est un produit de puissances positives
de variables d'amas qui appartiennent toutes au m\^eme amas.
La construction d'une \og base canonique\fg\ de l'alg\`ebre amass\'ee $\ca_Q$
est un probl\`eme important et encore tr\`es largement ouvert, voir
par exemple \cite{ShermanZelevinsky04} \cite{Dupont08} \cite{Cerulli09}.
On s'attend \`a ce qu'une telle base contienne tous les mon\^omes d'amas,
d'o\`u la conjecture~:

\begin{conj}[\cite{FominZelevinsky03}] \label{conj:independance}
Les mon\^omes d'amas sont lin\'eairement
ind\'ependants sur le corps $\Q$.
\end{conj}

Si $Q$ est un carquois de Dynkin, on sait \cite{CalderoKeller08} que
les mon\^omes d'amas forment une base de $\ca_Q$. Si $Q$ est {\em
acyclique}, c'est-\`a-dire n'admet aucun cycle orient\'e, la
conjecture r\'esulte d'un th\'eor\`eme de Geiss-Leclerc-Schr\"oer
\cite{GeissLeclercSchroeer07b}, qui montrent l'existence d'une \og base
g\'en\'erique\fg\ contenant les mon\^omes d'amas. La conjecture a aussi
\'et\'e d\'emontr\'ee pour des classes d'alg\`ebres amass\'ees \`a
coefficients (voir la
section~\ref{s:algebres-amassees-a-coefficients}), par exemple dans
les travaux \cite{FuKeller07} \cite{GeissLeclercSchroeer07b}
\cite{Demonet09}.

\begin{conj}[\cite{FominZelevinsky03}]  \label{conj:positivite}
Les variables d'amas s'\'ecrivent comme des polyn\^omes
de Laurent \`a coefficients entiers {\em positifs} en les variables
de tout amas.
\end{conj}

La cat\'egorification mono\"{\i}dale d\'evelopp\'ee par Leclerc
\cite{Leclerc08a} et Hernandez-Leclerc \cite{HernandezLeclerc09}
(voir la section~\ref{ss:categorification-monoidale}) a permis r\'ecemment de montrer cette
conjecture d'abord pour les carquois de type $A_n$ et $D_4$, voir
\cite{HernandezLeclerc09}, puis pour tout carquois admettant une
orientation bipartite \cite{Nakajima09}, c'est-\`a-dire une
orientation o\`u tout sommet est une source ou un puits. Elle est
d\'emontr\'ee de fa\c{c}on combinatoire par
Musiker-Schiffler-Williams \cite{MusikerSchifflerWilliams09} pour
tous les carquois associ\'es \`a des surfaces \`a bord marqu\'ees
\cite{FominShapiroThurston08} et par Di Francesco-Kedem \cite{DiFrancescoKedem09}
pour les carquois associ\'es au $T$-syst\`eme de type $A$.

Nous renvoyons \`a \cite{FominZelevinsky03a} et
\cite{FominZelevinsky07} pour de nombreuses autres conjectures sur
les alg\`ebres amass\'ees et \`a \cite{DerksenWeymanZelevinsky09}
pour la solution d'une bonne partie de ces conjectures gr\^ace \`a
la cat\'egorification (voir la section~\ref{ss:carquois-a-potentiel}).

\subsection{$Y$-graines, application \`a la conjecture de p\'eriodicit\'e}
\label{ss:Y-graines}

Soient $n\geq 1$ un entier et $\cg$ le corps $\Q(y_1, \ldots, y_n)$
engendr\'e par des ind\'etermi\-n\'ees~$y_i$. Une {\em $Y$-graine} est un couple $(R,v)$, o\`u
$R$ est un bon carquois et $v$ une suite \mbox{$v_1$, \ldots, $v_n$}
qui engendre librement le corps $\cg$ (nous nous \'ecartons quelque peu de la d\'efinition
dans \cite{FominZelevinsky07}). Si $(R,v)$ est une
$Y$-graine et $k$ un sommet de $R$, la {\em mutation $\mu_k(R,v)$}
est la $Y$-graine $(R',v')$, o\`u $R'=\mu_k(R)$ et
\[
v'_i=\left\{ \begin{array}{ll} v_i^{-1} & \mbox{si $i=k$,} \\
v_i (1+v_k)^m & \mbox{si le nombre de fl\`eches $i\to k$ est $m\geq 1$,} \\
v_i (1+v_k^{-1})^{-m} & \mbox{si le nombre de fl\`eches $k\to i$ est $m\geq 1$,}\\
v_i & \mbox{sinon.}
\end{array} \right.
\]
Par exemple, les $Y$-graines obtenues \`a partir de $y_1 \to y_2$ par
mutation sont, en \'ecrivant les variables $v_i$ \`a la place
des sommets $i$~:
\begin{align*}
\left( \frac{1}{y_1} \la \frac{y_1 y_2}{1+y_1} \right) &\maparr{\mu_2}
\left( \frac{y_2}{1+y_1+y_1y_2} \to \frac{1+y_1}{y_2} \right) \maparr{\mu_1}
\left( \frac{1+y_1+y_1y_2}{y_2} \la \frac{1}{y_1(1+y_2)} \right) \\
 &\maparr{\mu_2} \left( \frac{1}{y_2} \to y_1(1+y_2) \right) \maparr{\mu_1} \left( y_2 \la y_1 \right).
\end{align*}
Les $Y$-graines jouent un r\^ole important dans la th\'eorie de
Teichm\"uller sup\'erieure de Fock-Goncharov \cite{FockGoncharov03} \cite{FockGoncharov06a}
et dans l'\'etude par Kontsevich-Soibelman \cite{KontsevichSoibelman08}
des invariants de Donaldson-Thomas des carquois \`a potentiel \cite{DerksenWeymanZelevinsky08}.
Elles sont li\'ees aux $X$-graines par des conjectures de dualit\'e
\cite{FockGoncharov03} \'etudi\'ees syst\'ematiquement par Fomin-Zelevinsky
dans \cite{FominZelevinsky07}. En particulier, dans \cite{FominZelevinsky07},
les auteurs montrent que les
deux types de graines se d\'eterminent mutuellement si, en
m\^eme temps que $\ca_Q$, on consid\`ere aussi $\ca_{\tilde{Q}}$,
o\`u $\tilde{Q}$ est {\em l'extension principale} de $Q$ obtenue
\`a partir de $Q$ en rajoutant de nouveaux sommets $n+1$, \ldots, $2n$
et une nouvelle fl\`eche $(n+i) \to i$ pour tout $1\leq i\leq n$.
Ces liens combin\'es avec la cat\'egorification additive (voir
section~\ref{ss:categorification-additive-categorie-amassee})
ont permis r\'ecemment une application des alg\`ebres amass\'ees
\`a l'\'etude de syst\`emes dynamiques discrets issus de la
physique math\'ematique.

Soient $\Delta$ et $\Delta'$ deux
diagrammes de Dynkin simplement lac\'es. Notons $1, \ldots, n$
et $1, \ldots, n'$ leurs sommets et $A$ et $A'$ leurs matrices
d'incidence, le coefficient en position $(i,j)$ valant $1$ s'il
existe une ar\^ete entre $i$ et $j$ et $0$ sinon.
Notons $h$ et $h'$ les nombres de Coxeter de $\Delta$ et $\Delta'$.
Le {\em $Y$-syst\`eme} associ\'e \`a $\Delta$ et $\Delta'$ est un syst\`eme
infini d'\'equations de r\'ecurrence en des variables $Y_{i,j,t}$
associ\'ees aux sommets $(i,j)$ du produit $\Delta \times \Delta'$
et d\'ependant d'un param\`etre de temps discret $t\in\Z$. Les
\'equations du $Y$-syst\`eme sont
\begin{equation} \label{eq:Y-system}
Y_{i,i',t-1} Y_{i,i',t+1} =
\frac{\prod_{j=1}^n (1+Y_{j,i',t})^{a_{ij}}}{\prod_{j'=1}^{n'} (1+Y_{i,j',t}^{-1})^{a'_{i'j'}}} \ko
\end{equation}
pour tout sommet $(i,i')$ du produit et tout entier $t$.

\begin{theo} \label{thm:periodicite} Toutes les solutions du $Y$-syst\`eme sont p\'eriodiques
par rapport au param\`etre $t$ de p\'eriode divisant $2(h+h')$.
\end{theo}

Ce th\'eor\`eme vient confirmer la \og conjecture de p\'eriodicit\'e\fg\
formul\'ee par Al.~B.~Za\-mo\-lod\-chi\-kov \cite[(12)]{Zamolodchikov91}
pour $\Delta'=A_1$, par Kuniba-Nakanishi \cite[(2a)]{KunibaNakanishi92}
pour $\Delta'=A_m$ et par Ravanini-Valleriani-Tateo
\cite[(6.2)]{RavaniniTateoValleriani93} dans le cas g\'en\'eral. Le
th\'eor\`eme a \'et\'e d\'emontr\'e
\begin{itemize}

\item pour $(A_n, A_1)$  par Frenkel-Szenes \cite{FrenkelSzenes95}
(qui donnent des solutions explicites) et par Gliozzi-Tateo
\cite{GliozziTateo96} (\`a l'aide de calculs de volumes de
$3$-vari\'et\'es)~;

\item par Fomin-Zelevinsky \cite{FominZelevinsky03b} pour $(\Delta, A_1)$,
o\`u $\Delta$ n'est pas n\'ecessairement simplement lac\'e (ils utilisent
les m\'ethodes de leur th\'eorie des alg\`ebres amass\'ees et un
calcul sur ordinateur pour les types exceptionnels; ce calcul
peut \^etre \'evit\'e maintenant gr\^ace \`a \cite{YangZelevinsky09})~;

\item pour $(A_n, A_m)$ par Volkov \cite{Volkov07}, qui construit des
  solutions explicites gr\^ace \`a des consid\'erations
  de g\'eom\'etrie projective \'el\'ementaire, et par Szenes
  \cite{Szenes06}, qui interpr\`ete le syst\`eme comme un syst\`eme de
  connexions plates sur un graphe; la d\'emonstration d'un \'enonc\'e \'equivalent est due
  \`a Henriques \cite{Henriques07}~;

\item pour $(\Delta, \Delta')$ quelconques dans \cite{Keller08c} \cite{Keller10a} \`a l'aide de
la cat\'egorification additive, voir la
section~\ref{ss:categorification-additive-categorie-amassee}.

% \item par Hernandez-Leclerc pour $(A_n, A_1)$ \`a l'aide de repr\'esentations
% d'alg\`ebres quantiques affines (qui donnent des formules pour
% les solutions en termes de $q$-caract\`eres).

\end{itemize}

Pour des diagrammes non simplement lac\'es $\Delta$ et $\Delta'$, deux
variantes g\'en\'eralis\'ees de la conjecture existent~: la premi\`ere
se ram\`ene au th\'eor\`eme~\ref{thm:periodicite} par la technique du
\og pliage\fg, voir \cite{FominZelevinsky03b}~; la deuxi\`eme,
formul\'ee par Kuniba-Nakanishi \cite[(2a)]{KunibaNakanishi92} et
Kuniba-Nakanishi-Suzuki \cite[B.6]{KunibaNakanishiSuzuki94}, fait
intervenir le double de la somme des nombres de Coxeter {\em duaux}
(voir par exemple le chapitre~6 de \cite{Kac90})~; elle a
\'et\'e d\'emontr\'ee dans 
\cite{InoueIyamaKellerKunibaNakanishi10a}
\cite{InoueIyamaKellerKunibaNakanishi10b}.

\section{ALG\`EBRES AMASS\'EES \`A COEFFICIENTS}
\label{s:algebres-amassees-a-coefficients}

Nous allons g\'en\'eraliser l\'eg\`erement la d\'efinition
donn\'ee \`a la section~\ref{s:algebres-amassees-associees-aux-carquois}
pour obtenir la classe des \og alg\`ebres
amass\'ees antisym\'etriques de type g\'eom\'etrique\fg.
Cette classe contient de nombreuses alg\`ebres d'origine
g\'eom\'etrique munies de \og bases duales semi-canoniques\fg.
La construction d'une grande partie d'une telle base est l'une
des applications les plus remarquables des alg\`ebres amass\'ees.

Nous renvoyons \`a \cite{FominZelevinsky07}
pour la d\'efinition des \og alg\`ebres amass\'ees antisym\'etrisables
\`a coefficients dans un semi-corps\fg, qui constituent la classe
la plus g\'en\'erale consid\'er\'ee jusqu'\`a maintenant.

\subsection{D\'efinition} \label{ss:def-cluster-alg-coeff}
Soient $1 \leq n \leq m$ des entiers. Soit $\tilde{Q}$ un
{\em carquois glac\'e de type $(n,m)$}, c'est-\`a-dire
un bon carquois \`a $m$ sommets
et qui ne comporte aucune fl\`eche entre sommets $i,j$
tous les deux strictement plus grands que $n$.
La {\em partie principale} de $\tilde{Q}$ est le sous-carquois
plein $Q$ dont les sommets sont $1$, \ldots, $n$ (un sous-carquois
est {\em plein} si, avec deux sommets, il contient toutes
les fl\`eches qui les relient). Les sommets $n+1$, \ldots, $m$
sont les {\em sommets gel\'es}.
{\em L'alg\`ebre amass\'ee associ\'ee au carquois glac\'e $\tilde{Q}$}
\[
\ca_{\tilde{Q}} \subset \Q(x_1, \ldots, x_m)
\]
est d\'efinie de la m\^eme fa\c{c}on que l'alg\`ebre amass\'ee
associ\'ee \`a un carquois
(section~\ref{s:algebres-amassees-associees-aux-carquois})
sauf que
\begin{itemize}
\item seules les mutations par rapport \`a des sommets non gel\'es
sont admises et aucune fl\`eche entre sommets gel\'es
n'est introduite lors des mutations~;
\item les variables $x_{n+1}$, \ldots, $x_m$, qui font partie
de tous les amas, sont appel\'ees {\em coefficients} plut\^ot
que variables d'amas~;
\item le {\em type amass\'e} du carquois glac\'e est celui
de sa partie principale (s'il est d\'efini).
\end{itemize}
Souvent, on consid\`ere des localisations de $\ca_{\tilde{Q}}$
obtenues en inversant certains des coef\-ficients.
Si $K$ est une
extension de $\Q$ et $A$ une $K$-alg\`ebre (associative avec $1$),
{\em une structure d'alg\`ebre amass\'ee \`a coefficients de type
  $\tilde{Q}$ sur $A$} est la donn\'ee d'un isomorphisme $\phi$
de $\ca_{\tilde{Q}} \ten_\Q K$ sur $A$. Un tel isomorphisme
est d\'etermin\'e par les images des coefficients et des
variables de la graine initiale $\phi(x_i)$, $1 \leq i\leq m$. Nous appellerons
{\em graine initiale de $A$} la donn\'ee du carquois
$\tilde{Q}$ et des $\phi(x_i)$.

\subsection{Exemple~: le c\^one sur la grassmannienne des plans d'un espace vectoriel}

Soit $n\geq 1$ un entier. Soit $A$ l'alg\`ebre des fonctions
polynomiales sur le c\^one au-dessus de la grassmannienne des plans de
$\C^{n+3}$.  Cette alg\`ebre est engendr\'ee par les coordonn\'ees de
Pl\"ucker $x_{ij}$, $1\leq i<j\leq n+3$, assujetties aux relations de
Pl\"ucker~: pour tout quadruplet d'entiers $i<j<k<l$ compris entre
$1$ et $n+3$, nous avons
\begin{equation} \label{eq:Plucker}
x_{ik} x_{jl} = x_{ij} x_{kl} + x_{jk} x_{il}.
\end{equation}
% Tout plan $P$ de $\C^{n+3}$ donne une droite dans ce
% c\^one, \`a savoir celle engendr\'ee par la suite des
% mineurs $2\times 2$ d'une matrice $(n+3)\times 2$ dont
% les colonnes engendrent $P$.
Notons que les mon\^omes dans cette relation sont naturellement associ\'es
aux diagonales et aux c\^ot\'es du carr\'e
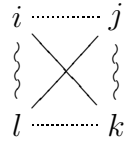
\[
\xymatrix{ i \ar@{.}[r] \ar@{~}[d] \ar@{-}[dr] & j \ar@{-}[dl] \ar@{~}[d] \\
l \ar@{.}[r] & k}
\]
L'id\'ee est d'interpr\'eter cette relation comme une
relation d'\'echange dans une alg\`ebre amass\'ee (\`a coefficients).
Pour d\'ecrire cette alg\`ebre, consid\'erons, dans le plan
affine euclidien, un polygone r\'egulier $P$
dont les sommets
sont num\'erot\'es de $1$ \`a $n+3$.
Consid\'erons la variable $x_{ij}$ comme associ\'ee au
segment $[ij]$ joignant les sommets $i$ et $j$.

\begin{proposition}[\protect{\cite[Example 12.6]{FominZelevinsky03}}]
\label{prop:plans} L'alg\`ebre $A$ a une structure d'alg\`ebre
amass\'ee \`a coefficients telle que
\begin{itemize}
\item[-] les coefficients soient les variables $x_{ij}$ associ\'ees aux
c\^ot\'es de $P$~;
\item[-] les variables d'amas soient les variables $x_{ij}$ associ\'ees
aux diagonales de $P$~;
\item[-] les amas soient les $n$-uplets de variables d'amas
correspondant \`a des diagonales qui forment une triangulation de $P$.
\end{itemize}
En outre, les relations d'\'echange sont exactement les relations
de Pl\"ucker et le type amass\'e est $A_n$.
\end{proposition}

Une triangulation de $P$ d\'etermine une graine initiale
pour l'alg\`ebre amass\'ee et les relations d'\'echange v\'erifi\'ees
par les variables d'amas initiales d\'eterminent le carquois glac\'e $\tilde{Q}$.
Par exemple, on v\'erifie que, dans le dessin suivant, la triangulation et le carquois glac\'e
(dont les sommets gel\'es sont entour\'es de bo\^{\i}tes) se correspondent
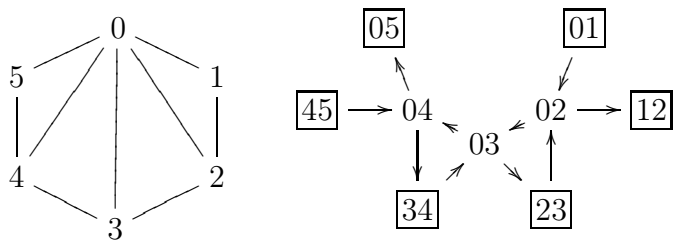
\[
\begin{xy} 0;<0.5pt,0pt>:<0pt,-0.5pt>::
(76,0) *+{0} ="0",
(150,37) *+{1} ="1",
(150,112) *+{2} ="2",
(74,151) *+{3} ="3",
(0,112) *+{4} ="4",
(0,37) *+{5} ="5",
(275,0) *+{\framebox[3ex]{05}} ="6",
(300,62) *+{04} ="7",
(350,87) *+{03} ="8",
(400,62) *+{02} ="9",
(425,0) *+{\framebox[3ex]{01}} ="10",
(225,62) *+{\framebox[3ex]{45}} ="11",
(300,137) *+{\framebox[3ex]{34}} ="12",
(400,137) *+{\framebox[3ex]{23}} ="13",
(475,62) *+{\framebox[3ex]{12}} ="14",
"0", {\ar@{-}"1"},
"0", {\ar@{-}"2"},
"0", {\ar@{-}"3"},
"0", {\ar@{-}"4"},
"5", {\ar@{-}"0"},
"1", {\ar@{-}"2"},
"2", {\ar@{-}"3"},
"3", {\ar@{-}"4"},
"4", {\ar@{-}"5"},
"7", {\ar"6"},
% "6", {\ar"11"},
"8", {\ar"7"},
"11", {\ar"7"},
"7", {\ar"12"},
"9", {\ar"8"},
"12", {\ar"8"},
"8", {\ar"13"},
"10", {\ar"9"},
"13", {\ar"9"},
"9", {\ar"14"},
% "14", {\ar"10"},
\end{xy}
\]
De nombreuses autres alg\`ebres de coordonn\'ees (homog\`enes) de
vari\'et\'es alg\'ebriques classiques admettent \'egalement des
structures d'alg\`ebres amass\'ees (sup\'erieures, voir la
section~\ref{ss:max-unipotent-sln}), notamment les grassmanniennes
\cite{Scott06} et les doubles cellules de Bruhat
\cite{BerensteinFominZelevinsky05}. Certaines de ces alg\`ebres
n'ont qu'un nombre fini de variables d'amas et donc un type amass\'e
bien d\'efini. Voici quelques exemples extraits de \cite{FominZelevinsky03a},
o\`u $N$ est un sous-groupe unipotent maximal~:

\begin{tabular}{|c|c|c|c|c|c|c|c|c|c|} \hline
$Gr_{2,n+3}$ & $Gr_{3,6}$ & $Gr_{3,7}$ & $Gr_{3,8}$ & $SL_3/N$ & $SL_4/N$ & $SL_5/N$ & $Sp_4/N$ & $SL_2$ & $SL_3$ \\ \hline
$A_n$ & $D_4$ & $E_6$ & $E_8$ & $A_1$ & $A_3$ & $D_6$ & $B_2$ & $A_1$ & $D_4$  \\ \hline
\end{tabular}
\bigskip

Un analogue de la proposition~\ref{prop:plans} pour les doubles cellules
de Bruhat r\'eduites \cite{BerensteinZelevinsky01} est d\^u \`a
Yang et Zelevinsky \cite{YangZelevinsky09}. Ils obtiennent ainsi une
alg\`ebre amass\'ee (\`a coefficients principaux) avec une description
explicite des variables d'amas pour tout diagramme de Dynkin.

\subsection{Exemple~: le sous-groupe unipotent maximal de $SL(n+1,\C)$}
\label{ss:max-unipotent-sln}

Soient $n$ un entier positif et $N$ le sous-groupe de $SL(n+1,\C)$
form\'e des matrices triangulaires sup\'erieures dont les coefficients
diagonaux sont tous \'egaux \`a $1$. Pour $1\leq i,j\leq n+1$ et
$g\in N$, soit $F_{ij}(g)$ la sous-matrice carr\'ee de taille
maximale de $g$ qui comporte le coefficient $g_{ij}$ dans son coin inf\'erieur gauche.
Soit $f_{ij}(g)$ le d\'eterminant de $F_{ij}(g)$. Nous
consid\'erons les fonctions polynomiales $f_{ij}: N \to \C$
pour $1\leq i\leq n$ et $i+j\leq n+2$.
{\em L'alg\`ebre amass\'ee sup\'erieure} associ\'ee \`a un
carquois glac\'e $\tilde{Q}$ \`a $m$ sommets est la sous-alg\`ebre de
$\Q(x_1, \ldots, x_m)$ form\'ee des \'el\'ements qui s'expriment
comme des polyn\^omes de Laurent en les variables de tout
amas associ\'e \`a $\tilde{Q}$.

\begin{theo}[\cite{BerensteinFominZelevinsky05}]
\label{thm:alg-amassee-matrices-unipotentes}
L'alg\`ebre des
fonctions polynomiales $\C[N]$ a une structure d'alg\`ebre amass\'ee
sup\'erieure dont la graine initiale est donn\'ee par
\end{theo}
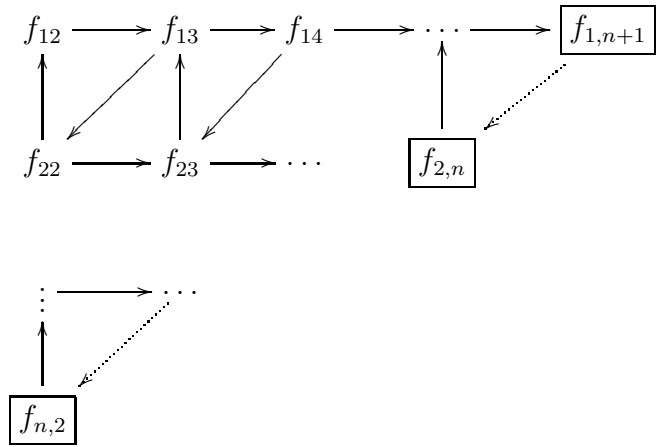
\[
\xymatrix{
f_{12} \ar[r] & f_{13} \ar[r] \ar[dl] & f_{14} \ar[dl] \ar[r] & \ldots \ar[r] & \framebox{$f_{1,n+1}$} \ar@{.>}[dl] \\
f_{22} \ar[u] \ar[r]
             & f_{23} \ar[u] \ar[r]  & \ldots         & \framebox{$f_{2,n}$} \ar[u]     &   \\
\vdots \ar[r]       &  \ldots \ar@{.>}[dl]                     &          &               &   \\
\framebox{$f_{n,2}$} \ar[u] &                    &                &               &
}
\]

Il n'est pas difficile de v\'erifier que cette structure est
de type amass\'e $A_3$ pour $n=3$, $D_6$ pour $n=4$ et qu'elle a une
infinit\'e de variables d'amas pour $n\geq 5$.

Un th\'eor\`eme de Fekete \cite{Fekete1912}, g\'en\'eralis\'e dans
\cite{BerensteinFominZelevinsky96}, affirme qu'une matrice carr\'ee \`a
$n+1$ lignes et \`a coefficients r\'eels est {\em totalement positive}
(\ie tous ses mineurs sont $>0$) si les $(n+1)^2$ mineurs suivants
sont strictement positifs~: tous les mineurs form\'es des 
$k$~premi\`eres lignes et de $k$ colonnes cons\'ecutives pour $1\leq
k\leq n+1$. Une matrice $g\in N$ \`a coefficients r\'eels est
{\em totalement positive} si tous les mineurs non identiquement nuls
sur $N$ sont strictement positifs en $g$. Cette condition
est \'equivalente \`a ce que l'on ait $f_{ij}(g)>0$ pour les
$f_{ij}$ de la graine initiale du th\'eor\`eme. Comme les
relations d'\'echange ne font pas intervenir de soustraction,
tout amas non initial $C$ donne \'egalement un crit\`ere de
positivit\'e~: la matrice $g\in N$ est totalement positive si et seulement si
l'on a $u_{ij}(g)>0$ pour toute variable d'amas $u_{ij}$ dans $C$.

\subsection{Exemple d'application aux bases duales semi-canoniques}
\label{ss:application-aux-bases-duales-semi-canoniques}
L'exemple~\ref{ss:max-unipotent-sln} se g\'en\'eralise. Soit, en effet, $G$
un groupe alg\'ebrique semi-simple complexe et $N$ un sous-groupe
unipotent maximal de $G$. Alors d'apr\`es \cite{BerensteinFominZelevinsky05},
l'alg\`ebre $\C[N]$ est munie d'une structure canonique
d'alg\`ebre amass\'ee sup\'erieure.
Soit $\mathfrak{n}$ l'alg\`ebre de Lie du groupe alg\'ebrique $N$.
Dans \cite{Lusztig00}, Lusztig construit une base distingu\'ee, la
{\em base semi-canonique}, de (l'espace vectoriel complexe sous-jacent \`a)
l'alg\`ebre enveloppante $U(\mathfrak{n})$.
Le dual restreint de la cog\`ebre $U(\mathfrak{n})$
est canoniquement isomorphe \`a $\C[N]$, qui est donc muni
de la base duale de celle construite par Lusztig, appel\'ee
{\em base duale semi-canonique}.

\begin{theo}[Geiss-Leclerc-Schr\"oer \cite{GeissLeclercSchroeer06}]
\label{thm:monomes-damas-base-semi-canonique}
  Tout mon\^ome d'amas de $\C[N]$ (munie de la structure d'alg\`ebre
  amass\'ee du
  th\'eor\`eme~\ref{thm:alg-amassee-matrices-unipotentes}) fait partie
  de la base duale semi-canonique.
\end{theo}

Notons que ce th\'eor\`eme implique la conjecture~\ref{conj:independance}
sur l'ind\'ependance des mon\^omes d'amas pour cette classe d'alg\`ebres amass\'ees.
L'alg\`ebre $U(\mathfrak{n})$ est \'egalement munie de
la {\em base canonique} obtenue par sp\'ecialisation \`a
partir de la base canonique \cite{Kashiwara90} \cite{Lusztig90}
du groupe quantique $U_q(\mathfrak{n})$.

\begin{theo}[Geiss-Leclerc-Schr\"oer \cite{GeissLeclercSchroeer05}]
Pour $G=SL(n+1,\C)$, la base cano\-nique
co\"{\i}ncide avec la base semi-canonique de $U(\mathfrak{n})$
si et seulement si $n\leq 4$.
\end{theo}

N\'eanmoins, Geiss-Leclerc-Schr\"oer conjecturent qu'au moins les
\og parties rigides\fg\ des bases duales canonique et semi-canonique
co\"{\i}ncident. Plus pr\'ecis\'ement, un cas parti\-culier
de la conjecture~23.2 de \cite{GeissLeclercSchroeer07b} nous donne
la conjecture qui suit.

\begin{conj}[Geiss-Leclerc-Schr\"oer \cite{GeissLeclercSchroeer07b}]
  Tout mon\^ome d'amas de $\C[N]$ appartient aussi \`a la base duale
  canonique.
\end{conj}

Dans un travail de longue haleine qui a abouti \`a
\cite{GeissLeclercSchroeer07b}, Geiss-Leclerc-Schr\"oer ont
g\'en\'eralis\'e les
th\'eor\`emes~\ref{thm:alg-amassee-matrices-unipotentes} et
\ref{thm:monomes-damas-base-semi-canonique} de l'alg\`ebre $\C[N]$ aux
alg\`ebres de coordonn\'ees de cellules unipotentes de groupes de
Kac-Moody simplement lac\'es.  Nous renvoyons \`a
\cite{GeissLeclercSchroeer08a} pour une introduction et une synth\`ese
des r\'esultats dans le cas fini, et \`a \cite{Demonet09} pour une
extension (partielle) au cas non simplement lac\'e.

\section{CAT\'EGORIFICATIONS}
\label{s:categorifications}

Les d\'emonstrations du th\'eor\`eme de p\'eriodicit\'e
(Th\'eor\`eme~\ref{thm:periodicite}) et du th\'eor\`eme sur la base
duale semi-canonique
(Th\'eor\`eme~\ref{thm:monomes-damas-base-semi-canonique}) s'appuient
sur la \og cat\'egorification additive\fg\ des alg\`ebres amass\'ees~; celle
des cas connus de la conjecture de positivit\'e~\ref{conj:positivite}
sur la \og cat\'egorification mono\"{\i}dale\fg. Nous allons esquisser les
id\'ees principales de ces m\'ethodes. La
section~\ref{ss:carquois-a-potentiel} est consacr\'ee \`a la m\'ethode de
la cat\'egorification \`a l'aide des \og carquois \`a potentiel\fg. \`A ce
jour, c'est la seule m\'ethode qui permette de traiter des alg\`ebres
amass\'ees associ\'ees \`a des carquois finis arbitraires (sans
boucles ni $2$-cycles).

\subsection{Cat\'egorification additive~: la cat\'egorie amass\'ee}
\label{ss:categorification-additive-categorie-amassee}

Soit $Q$ un carquois fini d'ensemble de sommets $\{1, \ldots, n\}$.
Un {\em chemin de $Q$} est une composition formelle $(j|\alpha_s| \ldots | \alpha_1|i)$
d'un nombre $s$ positif ou nul de fl\`eches $\alpha_i$ telle que 
$s(\alpha_i)=t(\alpha_{i-1})$ pour $1\leq i\leq s$.
En particulier, pour tout sommet $i$, nous avons le {\em chemin paresseux
$e_i=(i|i)$} de longueur nulle, neutre pour la composition naturelle des chemins.
Une {\em repr\'esentation} (complexe) de $Q$ est la
donn\'ee $V$ d'espaces vectoriels complexes de dimension finie $V_i$,
$i\in Q_0$, et d'applications lin\'eaires $V_\alpha: V_i \to V_j$ pour
toute fl\`eche $\alpha: i\to j$ de $Q$. Une repr\'esentation est donc
un diagramme d'espaces vectoriels de la forme donn\'ee par $Q$. Un
{\em morphisme de repr\'esentations} est un morphisme de diagrammes.
On obtient ainsi la {\em cat\'egorie $\rep(Q)$} des repr\'esentations
de $Q$. C'est une cat\'egorie ab\'elienne \'equivalente
\`a la cat\'egorie des modules de $\C$-dimension finie sur une
alg\`ebre, \`a savoir {\em l'alg\`ebre des chemins $\C Q$} (une base
de cette alg\`ebre est form\'ee des chemins de $Q$~; le produit de
deux chemins composables est leur composition,
le produit de chemins non composables est nul). En particulier,
nous avons des notions naturelles de sous-repr\'esentation, de
repr\'esentation simple, de somme directe et de repr\'esentation
ind\'ecomposable ($=$ repr\'esentation non nulle qui n'est pas somme
directe de deux sous-repr\'esentations non nulles).

Supposons que $Q$ est un carquois de Dynkin de type $\Delta$. Alors d'apr\`es le
th\'eor\`eme de Gabriel \cite{Gabriel73}, on a une bijection de
l'ensemble des classes d'isomorphisme de repr\'esentations
ind\'ecomposables de $Q$ sur l'ensemble des racines
positives de $\Delta$~; \`a une repr\'esentation ind\'ecomposable
$V$, cette bijection associe la racine
$
\sum_{i=1}^n (\dim V_i) \alpha_i
$,
o\`u les $\alpha_i$ sont les racines simples. En composant cette bijection
avec celle de la partie~c) du th\'eor\`eme~\ref{thm:class-finie-detaillee}
de Fomin-Zelevinsky,
nous obtenons une bijection de l'ensemble des classes
d'isomorphisme de repr\'esentations ind\'ecomposables sur l'ensemble
des variables d'amas non initiales de l'alg\`ebre amass\'ee $\ca_Q$~:
\`a une repr\'esentation ind\'ecomposable $V$, cette bijection associe l'unique
variable d'amas non initiale $X_V$ dont le d\'enominateur est
$x_1^{d_1} \cdots x_n^{d_n}$, o\`u $d_i=\dim V_i$. Il est
remarquable que le num\'erateur de $X_V$ admette aussi une interpr\'etation
naturelle en termes de la repr\'esentation $V$.
Pour expliciter cette interpr\'etation, nous avons besoin
de quelques notations suppl\'ementaires~: soient $V$ une repr\'esentation
quelconque de $Q$ et $d_i=\dim V_i$, $i\in Q_0$.
Pour un \'el\'ement $e\in \N^n$, notons $Gr_e(V)$ l'ensemble
des sous-representations $U$ de $V$ telles que $\dim U_i = e_i$.
La donn\'ee d'un point de $Gr_e(V)$ est donc la donn\'ee d'une
famille de sous-espaces vectoriels $U_i \subset V_i$ telle que $U_i$
soit de dimension $e_i$
et que $V_\alpha(U_i)\subset U_j$ pour toute fl\`eche $\alpha: i\to j$
de $Q$. Cette description montre que $Gr_e(V)$ est une sous-vari\'et\'e
ferm\'ee du produit des grassmanniennes $Gr_{e_i}(V_i)$. En particulier,
c'est une vari\'et\'e projective (singuli\`ere en g\'en\'eral). Elle est appel\'ee
{\em grassmannienne des sous-repr\'esentations} (quiver Grassmannian) et \'etudi\'ee
dans \cite{CalderoReineke08}, par exemple. On note $\chi(Gr_e(V_i))$ la
caract\'eristique d'Euler-Poincar\'e de son espace topologique sous-jacent.
Posons
\[
CC(V)= \frac{1}{x_1^{d_1} \cdots x_n^{d_n}}
\sum_e \chi(Gr_e(V)) \prod_{i=1}^n x_i^{\sum_{j\to i} e_j + \sum_{i\to j} (d_j-e_j)} \ko
\]
o\`u les sommes dans l'exposant portent sur les fl\`eches de but~$i$
respectivement de source~$i$.
\begin{theo}[Caldero-Chapoton \cite{CalderoChapoton06}] Si $V$ est ind\'ecomposable,
nous avons $X_V = CC(V)$.
\end{theo}

Notons que la formule pour $CC(V)$ a un sens pour toute repr\'esentation
de tout carquois fini $Q$. Supposons maintenant que $Q$ est un carquois
sans cycles orient\'es quelconque. Une repr\'esentation $V$ de $Q$ est
{\em rigide} si son groupe d'auto-extensions $\Ext^1(V,V)$ dans la
cat\'egorie $\rep(Q)$ s'annule. La partie~c) du
th\'eor\`eme~\ref{thm:class-finie-detaillee} ne s'applique plus,
mais nous avons n\'eanmoins une param\'etrisation des variables
d'amas non initiales en termes de repr\'esentations de $Q$.

\begin{theo}[\cite{CalderoKeller06}] L'application $V \mapsto CC(V)$ induit une bijection de
l'ensemble des classes d'isomorphisme de repr\'esentations rigides
ind\'ecomposables de $Q$ sur l'ensemble des variables d'amas non initiales
de $\ca_Q$.
\end{theo}

Ce th\'eor\`eme fournit une interpr\'etation cat\'egorique de la
quasi-totalit\'e des variables d'amas de l'alg\`ebre
amass\'ee. Il se pose la question d'\'etendre cette interpr\'etation
aux variables initiales, aux relations d'\'echange et aux amas.
Pour ce faire, on agrandit la cat\'egorie des repr\'esentations~: soit $\cd_Q$ la
cat\'egorie d\'eriv\'ee born\'ee de la cat\'egorie ab\'elienne
$\rep(Q)$. Les objets de $\cd_Q$ sont donc les complexes born\'es de
repr\'esentations et ses morphismes sont obtenus \`a partir des
morphismes de complexes en inversant formellement les
quasi-isomorphismes. La cat\'egorie $\cd_Q$ est une cat\'egorie
triangul\'ee~; on note $\Sigma$ son foncteur suspension (qui n'est
autre que le foncteur de d\'ecalage des complexes $X\mapsto X[1]$).
Les ensembles de morphismes de $\cd_Q$ sont des espaces vectoriels de
dimension finie et $\cd_Q$ admet un {\em foncteur de Serre},
c'est-\`a-dire une auto-\'equivalence $S:\cd_Q \to \cd_Q$ telle qu'on
ait des isomorphismes bifonctoriels
\[
D\Hom(X,Y) = \Hom(Y,SX) \ko
\]
o\`u $D=\Hom_\C(?,\C)$ est la dualit\'e des espaces vectoriels complexes.
La {\em cat\'egorie amass\'ee} est la cat\'egorie d'orbites
\[
\cc_Q = \cd_Q/ (S^{-1} \circ \Sigma^2)^\Z
\]
de $\cd_Q$ sous l'action du groupe cyclique engendr\'e par l'automorphisme
$S^{-1}\circ \Sigma^2$. Elle est due \`a Buan-Marsh-Reineke-Reiten-Todorov
\cite{BuanMarshReinekeReitenTodorov06} et, de fa\c{c}on ind\'ependante
et sous une forme tr\`es diff\'erente, \`a
Caldero-Chapoton-Schiffler \cite{CalderoChapotonSchiffler06} pour
les carquois de Dynkin de type $A$.
La cat\'egorie $\cc_Q$ est canoniquement triangul\'ee \cite{Keller05}.
Ses espaces de morphismes sont de dimension finie et son foncteur
de Serre (induit par $S$) est isomorphe au carr\'e de son foncteur
suspension (induit par $\Sigma$). Cela signifie que $\cc_Q$
est {\em Calabi-Yau de dimension $2$}. Notons $\pi: \cd_Q \to \cc_Q$
le foncteur de projection canonique et
\[
\Ext^1(L,M)=\Hom(L,\Sigma M)
\]
pour des objets $L$ et $M$ de $\cc_Q$. Gr\^ace \`a la propri\'et\'e
de Calabi-Yau, on a $D\Ext^1(L,M)=\Ext^1(M,L)$. On appelle
{\em rigide} un objet $L$ tel que $\Ext^1(L,L)$ s'annule. Notons $P_i$ la
repr\'esentation qui correspond au $\C Q$-module $\C Q e_i$, $i\in Q_0$.
On peut montrer \cite{BuanMarshReinekeReitenTodorov06} que tout objet
$L$ de $\cc_Q$ se d\'ecompose de fa\c{c}on unique (\`a isomorphisme pr\`es)
sous la forme
\[
L=\pi(M) \oplus \bigoplus_{i\in Q_0} \Sigma \pi(P_i)^{m_i} \ko
\]
pour une repr\'esentation $M$ et des multiplicit\'es $m_i$. On pose
\[
CC(L) = CC(M) \cdot \prod_{i\in Q_0} x_i^{m_i}.
\]

\begin{theo}[\cite{CalderoKeller06}] \label{thm:categorie-amassee-bijection}
\begin{itemize}
\item[a)] On a $CC(L\oplus M) = CC(L)\cdot CC(M)$ pour tous $L$ et $M$ dans $\cc_Q$,
\item[b)] si $L$ et $M$ sont des objets de $\cc_Q$ tels que $\Ext^1(L,M)$ soit de
dimension $1$ et
\[
L \to E \to M \to \Sigma L \mbox{ et } M \to E' \to L \to \Sigma M
\]
soient deux triangles non scind\'es, on a
\begin{equation} \label{eq:formule-de-multiplication}
CC(L) \cdot CC(M) = CC(E) + CC(E').
\end{equation}
\item[c)]
L'application $CC$ induit une bijection de l'ensemble des objets rigides de $\cc_Q$
sur l'ensemble des mon\^omes d'amas de $\ca_Q$.
\item[d)] Par cette bijection, les
objets rigides ind\'ecomposables correspondent aux variables d'amas, et
un ensemble d'ind\'ecomposables rigides $T_1$, \ldots, $T_n$
correspond \`a un amas si et seulement si $\Ext^1(T_i,T_j)=0$ pour tous $i$, $j$.
\end{itemize}
\end{theo}

Les propri\'et\'es a) et b) fournissent une interpr\'etation
des relations d'\'echange.
La d\'emonstration du th\'eor\`eme est fond\'ee sur le travail de
plusieurs groupes d'auteurs~:
Buan-Marsh-Reiten-Todorov \cite{BuanMarshReitenTodorov07},
Buan-Marsh-Reiten \cite{BuanMarshReiten08},
Buan-Marsh-Reineke-Reiten-Todorov \cite{BuanMarshReinekeReitenTodorov06},
Marsh-Reineke-Zelevinsky \cite{MarshReinekeZelevinsky03},
\ldots\ et surtout Caldero-Chapoton \cite{CalderoChapoton06}.
Une autre d\'emonstration  de la
formule de multiplication~\ref{eq:formule-de-multiplication} est due
\`a Hubery \cite{Hubery06} pour des carquois dont le graphe
sous-jacent est un diagramme de Dynkin \'etendu, et
\`a Xiao-Xu \cite{XiaoXu07} dans le cas g\'en\'eral.

La construction de la cat\'egorie amass\'ee a \'et\'e
g\'en\'eralis\'ee des alg\`ebres $\C Q$ \`a une classe d'alg\`ebres de
dimension globale $2$ par Amiot \cite{Amiot08a}. Une version
g\'en\'eralis\'ee de l'application de Caldero-Chapoton et de la
formule de multiplication~\ref{eq:formule-de-multiplication} est due
\`a Palu \cite{Palu08a}. L'extension des r\'esultats de Palu au cas
de certaines alg\`ebres amass\'ees \`a coefficients est obtenue dans
\cite{FuKeller07}. La mutation dans une cat\'egorie $2$-Calabi-Yau
g\'en\'erale est construite par Iyama-Yoshino \cite{IyamaYoshino08}.
La d\'emonstration du th\'eor\`eme de
p\'eriodicit\'e~\ref{thm:periodicite} dans \cite{Keller08c} \cite{Keller10a} est
fond\'ee sur ces travaux.

\subsection{Cat\'egorification additive~: modules sur les alg\`ebres pr\'eprojectives}
\label{ss:categorification-additive-preprojective}

Nous allons d\'ecrire l'id\'ee de base de la d\'emonstration du
th\'eor\`eme~\ref{thm:monomes-damas-base-semi-canonique}.
Soient $\Delta$ un diagramme de Dynkin, $\mathfrak{g}$ l'alg\`ebre
de Lie simple complexe qui lui correspond et $\mathfrak{n}$
une sous-alg\`ebre nilpotente maximale de $\mathfrak{g}$.
Soit $N$ le groupe alg\'ebrique unipotent associ\'e \`a~$\mathfrak{n}$.
L'alg\`ebre de coordonn\'ees $\C[N]$ est le dual restreint
de la cog\`ebre $U(\mathfrak{n})$.
La {\em base duale semi-canonique} de l'alg\`ebre $\C[N]$ est
duale de la base semi-canonique de $U(\mathfrak{n})$ construite
par Lusztig \cite{Lusztig00}. Dans un premier temps,
nous allons d\'ecrire (suivant \cite{GeissLeclercSchroeer06})
la base duale semi-canonique en termes de
modules sur l'alg\`ebre pr\'eprojective~:
soit $Q$ un carquois de Dynkin de type $\Delta$. Soit $\ol{Q}$
le {\em double carquois}, obtenu \`a partir de $Q$ en rajoutant une fl\`eche
$\alpha^*: j\to i$ pour chaque fl\`eche $\alpha: i\to j$. Soit $\Lambda$
{\em l'alg\`ebre pr\'eprojective de~$Q$}, c'est-\`a-dire le quotient de
l'alg\`ebre des chemins $\C \ol{Q}$ (voir
la section~\ref{ss:categorification-additive-categorie-amassee})
par l'id\'eal bilat\`ere engendr\'e par la somme $\sum [\alpha, \alpha^*]$
prise sur l'ensemble des fl\`eches de $Q$. C'est une alg\`ebre de dimension
finie sur $\C$ qui est auto-injective (c'est-\`a-dire injective comme module sur elle-m\^eme).
Appelons {\em $\Lambda$-module} un $\Lambda$-module \`a gauche de dimension
finie sur $\C$. Pour un tel module $M$, son {\em vecteur dimension}
est la suite des entiers $\dim e_i M$, $i\in Q_0$. Soit $d$ une famille
d'entiers positifs index\'es par $Q_0$. On note
$\rep(\Lambda, d)$ la vari\'et\'e form\'ee des familles
de matrices
\[
M_\alpha : \C^{d_{s(\alpha)}} \to \C^{d_{t(\alpha)}} \ko \alpha\in \ol{Q}_1 \ko
\]
qui v\'erifient les relations de $\Lambda$, c'est-\`a-dire
d\'efinissent une structure de $\Lambda$-module sur la somme directe
des $\C^{d_i}$, $i\in Q_0$.  La vari\'et\'e $\rep(\Lambda, d)$ porte
une action naturelle par \og changement de base\fg\ du groupe $G_d = \prod
GL(d_i,\C)$ et les orbites de cette action sont en bijection avec les
classes d'isomorphisme de $\Lambda$-modules de vecteur dimension~$d$.
Notons $\cm_d$ l'espace vectoriel des fonctions constructibles et
$G_d$-invariantes sur la vari\'et\'e $\rep(\Lambda,d)$.  Notons
$U(\mathfrak{n})_d$ la composante gradu\'ee de $U(\mathfrak{n})$
associ\'ee au vecteur $\sum d_i \alpha_i$, o\`u les $\alpha_i$ sont
les racines simples.  Lusztig \cite{Lusztig00} a d\'efini une
injection lin\'eaire $\lambda_d$ de $U(\mathfrak{n})_d$ dans $\cm_d$.
Chaque $\Lambda$-module $M$ d\'efinit une forme lin\'eaire sur
$\cm_d$, \`a savoir la forme qui, \`a une fonction $f$, associe sa
valeur $f(M)$ en l'orbite d\'etermin\'ee par $M$. Par composition, le
module $M$ nous donne une forme lin\'eaire
\[
\delta_M: U(\mathfrak{n})_d \arr{\lambda_d} \cm_d \to \C.
\]
Or nous avons l'isomorphisme canonique $\iota : U(\mathfrak{n})_d^*
\to \C[N]_d$. Comme dans \cite{GeissLeclercSchroeer06}, nous posons
$\phi_M = \iota(\delta_M)$.  Nous obtenons ainsi une application
$M\mapsto \phi_M$ de la classe des $\Lambda$-modules dans l'alg\`ebre
$\C[N]$. On peut expliciter cette application en termes de
caract\'eristiques d'Euler-Poincar\'e de vari\'et\'es de drapeaux de
sous-repr\'esentations de~$M$, voir \cite{GeissLeclercSchroeer06}.
Cette description montre que l'application $M \mapsto \phi_M$ est
constructible sur la vari\'et\'e alg\'ebrique $\rep(\Lambda, d)$.
Donc chaque composante irr\'eductible de cette vari\'et\'e contient un
ouvert dense o\`u la fonction $M \mapsto \phi_M$ est
constante. On appelle {\em g\'en\'eriques} les modules $M$ appartenant
\`a de tels ouverts. Alors la base duale semi-canonique n'est autre que
\[
\{ \phi_M\; | \; \mbox{$M$ est g\'en\'erique} \}.
\]
Un module $M$ est {\em rigide} si l'espace $\Ext^1(M,M)$ s'annule. De fa\c{c}on
\'equivalente \cite{GeissLeclercSchroeer06}, l'orbite de $M$ dans $\rep(\Lambda,d)$,
o\`u $d$ est le vecteur dimension de $M$, est ouverte. En particulier, si $M$ est
rigide, alors il est g\'en\'erique et la fonction $\phi_M$ appartient \`a la base
duale semi-canonique. Pour montrer le th\'eor\`eme~\ref{thm:monomes-damas-base-semi-canonique},
il suffit donc de montrer que chaque mon\^ome d'amas est de la
forme $\phi_M$ pour un module rigide $M$. Pour cela, on proc\`ede par
r\'ecurrence~: on montre \cite{GeissLeclercSchroeer08b}
que les \'el\'ements de la graine initiale sont
des images de modules rigides ind\'ecomposables canoniques $T_1^{(0)}$,
$T_2^{(0)}$, \ldots, $T_m^{(0)}$. Puis on rel\`eve l'op\'eration de mutation
des amas \`a une classe convenable de suites $T_1$, \ldots, $T_m$ de modules
rigides ind\'ecomposables. Appelons {\em accessibles} les modules
rigides dont les facteurs directs ind\'ecomposables sont obtenus
par mutation it\'er\'ee \`a partir de la suite $T_1^{(0)}$,
$T_2^{(0)}$, \ldots, $T_m^{(0)}$. Le th\'eor\`eme suivant est
l'analogue pr\'ecis du th\'eor\`eme~\ref{thm:categorie-amassee-bijection}.
Ses parties a) et b) permettent
de relier la mutation des amas \`a la mutation des
modules rigides et donc d'effectuer la r\'ecurrence qui termine
la d\'emonstration du th\'eor\`eme~\ref{thm:monomes-damas-base-semi-canonique}.

\begin{theo}[Geiss-Leclerc-Schr\"oer \cite{GeissLeclercSchroeer06}] On a
\begin{itemize}
\item[a)] $\phi_{L\oplus M} = \phi_L \phi_M$,
\item[b)] Si $\Ext^1(L,M)$ est de dimension $1$ et que l'on a les
suites exactes non scind\'ees
\[
0 \to L \to E \to M \to 0 \mbox{ et } 0 \to M \to E' \to L \to 0 \ko
\]
alors on a $\phi_L \phi_M = \phi_E + \phi_{E'}$.
\item[c)] L'application $M \mapsto \phi_M$ induit
une bijection de l'ensemble des classes d'isomorphisme de modules rigides
accessibles sur l'ensemble des mon\^omes d'amas.
\item[d)] Les rigides ind\'ecomposables accessibles correspondent
aux variables d'amas et une suite $T_1$, \ldots, $T_m$ de tels
modules correspond \`a un amas si et seulement si l'on a
$\Ext^1(T_i, T_j)=0$ pour tous $i$ et $j$.
\end{itemize}
\end{theo}

\subsection{Cat\'egorification mono\"{\i}dale}
\label{ss:categorification-monoidale}
Les cat\'egorifications additives d\'ecrites ci-dessus se sont av\'er\'ees
tr\`es utiles et on sait les construire pour de grandes classes d'alg\`ebres
amass\'ees. De l'autre c\^ot\'e, elles semblent difficiles \`a exploiter
pour d\'emontrer la conjecture de positivit\'e~\ref{conj:positivite}, et la notion m\^eme de
cat\'egorification additive semble peu naturelle. La cat\'egorification
mono\"{\i}dale, introduite par Leclerc \cite{Leclerc08a} et
Hernandez-Leclerc \cite{HernandezLeclerc09}, consiste \`a r\'ealiser
une alg\`ebre amass\'ee comme l'anneau de Grothendieck d'une cat\'egorie
ab\'elienne mono\"{\i}dale. Elle est donc tr\`es naturelle. En outre,
comme nous allons le voir, son existence donne imm\'ediatement la conjecture
de positivit\'e~\ref{conj:positivite} (et la conjecture d'ind\'ependance~\ref{conj:independance}).
De l'autre c\^ot\'e, les cat\'egorifications mono\"{\i}dales
semblent tr\`es difficiles \`a construire.

Les notions suivantes \cite{Leclerc08a} sont fondamentales pour la suite~:
un objet simple $S$ d'une cat\'egorie ab\'elienne mono\"{\i}dale est {\em premier}
s'il n'admet pas de factorisation tensorielle non triviale~; il est {\em r\'eel}
si son carr\'e tensoriel est encore simple.

Soient $\ca$ une alg\`ebre amass\'ee \`a coefficients et $\ca_\Z$ son
sous-anneau engendr\'e par les variables d'amas et les coefficients.
Suivant \cite{HernandezLeclerc09},
une {\em cat\'egorification mono\"{\i}dale} de $\ca$ est la donn\'ee
d'une cat\'egorie ab\'elienne mono\"{\i}dale $\cm$ et d'un isomorphisme d'anneaux
\[
\phi: \ca_\Z \iso K_0(\cm)
\]
tel que $\phi$ induise
\begin{itemize}
\item[a)] une bijection de l'ensemble des mon\^omes d'amas sur l'ensemble
des classes d'objets simples r\'eels de $\cm$ et
\item[b)] une bijection de l'ensemble des variables d'amas et des coefficients
sur l'ensemble des classes d'objets simples, r\'eels et premiers de $\cm$.
\end{itemize}
Le tableau suivant, extrait de \cite{Leclerc08a}, r\'esume les correspondances
entre les structures associ\'ees \`a une alg\`ebre amass\'ee et
leurs rel\`evements dans une cat\'egorification additive respectivement mono\"{\i}dale.
%\begin{table}
\begin{center}
\begin{tabular}{|c|c|c|} \hline
alg\`ebre amass\'ee $\ca$ & cat\'egorification additive $\cc$
                             & cat\'egorification mono\"{\i}dale $\cm$ \\ \hline
$+$ & ? & $\oplus$ \\
$\times$ & $\oplus$ & $\otimes$ \\
mon\^ome d'amas & objet rigide & objet simple r\'eel \\
variable d'amas & ind\'ecomposable rigide & simple premier r\'eel \\ \hline
\end{tabular}
\bigskip
%\caption{Correspondences between categorifications}
%\label{table:correspondences}
%\end{table}
\end{center}
L'existence d'une cat\'egorification mono\"{\i}dale $\cm$ d'une alg\`ebre
amass\'ee $\ca$ a des cons\'equences tr\`es fortes pour $\ca$~: en effet,
l'alg\`ebre $\ca$ est alors munie d'une \og base canonique\fg, \`a savoir la base fournie par
les objets simples de $\cm$ et cette base contient les mon\^omes d'amas
car ceux-ci correspondent bijectivement aux classes dans $K_0(\cm)$ de certains
objets simples. En particulier,
la conjecture d'ind\'ependance est v\'erifi\'ee pour $\ca$.
De m\^eme, la conjecture de positivit\'e est v\'erifi\'ee pour $\ca$. En effet,
si on exprime une variable d'amas $x$ comme polyn\^ome de Laurent
\[
x=\frac{P(u_1, \ldots, u_m)}{u_1^{d_1} \cdots u_m^{d_m}}
\]
en les variables d'un amas $u_1$, \ldots, $u_m$, alors les coefficients
de $P$ sont les multiplicit\'es de certains objets simples dans
la classe du produit tensoriel $\phi(x u_1^{d_1} \cdots u_m^{d_m})$ et sont
donc des entiers positifs.

L'existence d'une cat\'egorification mono\"{\i}dale $\phi: \ca_\Z \to K_0(\cm)$
donne \'egalement des renseignements pr\'ecieux sur la structure mono\"idale
de $\cm$~: en effet, elle montre que le comportement des objets simples r\'eels
de $\cm$ est gouvern\'e par la combinatoire des amas de $\ca$.

Dans \cite{HernandezLeclerc09}, Hernandez-Leclerc exhibent des
cat\'egorifications mono\"{\i}dales conjecturales $\cm_l$ pour les
alg\`ebres amass\'ees $\ca_l$ associ\'ees \`a certains carquois
glac\'es $\tilde{Q}(\Delta,l)$, o\`u $\Delta$ est un diagramme de Dynkin
simplement lac\'e et $l\in\N$ un \og niveau\fg. Voici l'exemple du
carquois $\tilde{Q}(D_5, 3)$, o\`u les sommets gel\'es sont
marqu\'es par des $\bullet$.
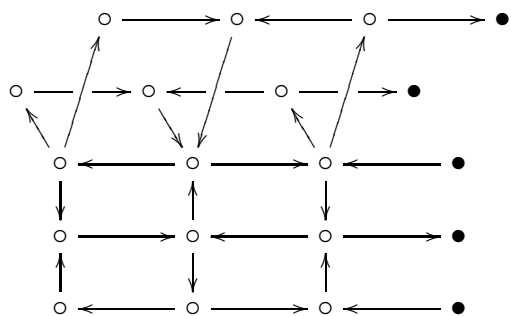
\[
\xymatrix@C=0.1cm@R=0.5cm{ & & \circ \ar[rrr] &  &  & \circ \ar[ldd] & & &
\circ \ar[lll] \ar[rrr] & & & \bullet  \\
     \circ \ar[rrr]|!{[dr];[urr]}\hole & &  & \circ \ar[rd] & & &
\circ \ar[lll]|!{[lu];[dll]}\hole \ar[rrr]|!{[rd];[rru]}\hole & & & \bullet  & &  \\
           & \circ \ar[lu] \ar[ruu] \ar[d] & & & \circ \ar[lll] \ar[rrr] & & &
\circ \ar[ruu] \ar[lu] \ar[d] & & & \bullet \ar[lll] & \\
           & \circ \ar[rrr] & & & \circ \ar[d] \ar[u] & & &
\circ \ar[lll] \ar[rrr] & & & \bullet  & \\
           & \circ \ar[u] & & & \circ \ar[lll] \ar[rrr] &  & &
\circ \ar[u] & & & \bullet \ar[lll] &
}
\]
Hernandez-Leclerc construisent les cat\'egories $\cm_l$ comme des sous-cat\'egories
mono\"{\i}dales de la cat\'egorie des repr\'esentations de dimension finie de
l'alg\`ebre affine quantique $U_q(\hat{\mathfrak{g}})$ associ\'ee \`a $\Delta$.
Ils construisent un morphisme d'anneaux $\phi: (\ca_l)_\Z \to K_0(\cm_l)$ qui
envoie les variables de l'amas initial sur les classes de certains modules
de Kirillov-Reshetikhin et conjecturent que $\phi$ est une cat\'egorification
mono\"{\i}dale de $\ca_l$ (Conjecture~13.2 de \cite{HernandezLeclerc09}).
Ils d\'emontrent leur conjecture pour $l\leq 1$ et $\Delta$ de
type $A_n$, $n\geq 1$, ou $D_4$ ainsi que pour $\Delta$ de type $A_2$ et $l=2$
(et observent que pour $\Delta=A_1$ et tout $l\in \N$, la conjecture r\'esulte du
travail de Chari-Pressley \cite{ChariPressley91}).

Dans \cite{Nakajima09}, Nakajima construit
des cat\'egorifications mono\"{\i}dales conjecturales $\cn_l$ pour les carquois
$\tilde{Q}(R,l)$ associ\'es \`a un carquois $R$ bipartite
et un niveau $l\in\N$. Les cat\'egories $\cn_l$ sont r\'ealis\'ees comme des
cat\'egories de faisceaux pervers sur des vari\'et\'es de carquois gradu\'es
\cite{Nakajima01} munies du produit tensoriel construit g\'eom\'etriquement dans \cite{VaragnoloVasserot03}.
Si $R$ est un carquois de Dynkin et $l=1$, la cat\'egorie $\cn_l$ est \'equivalente
\`a $\cm_l$ et Nakajima d\'emontre qu'elle est une cat\'egorification
mono\"{\i}dale de l'alg\`ebre $\ca_l$ confirmant ainsi la conjecture
de Hernandez-Leclerc. Pour un carquois bipartite $R$ quelconque et
$l=1$, il montre que l'anneau $(\ca_l)_\Z$ se plonge dans $K_0(\cn_l)$
de telle fa\c{c}on que les mon\^omes d'amas sont envoy\'es sur
des objets simples. Ceci entra\^{\i}ne la conjecture de positivit\'e pour
$\ca_l$.

\subsection{Cat\'egorification via les carquois \`a potentiel} \label{ss:carquois-a-potentiel}

Inspir\'es par des travaux de physiciens (voir par exemple
la section~6 dans \cite{FengHananyHeUranga01}) Derksen-Weyman-Zelevinsky
ont \'etendu \cite{DerksenWeymanZelevinsky08} l'op\'eration
de mutation des carquois aux carquois \`a potentiel
et leurs repr\'esentations d\'ecor\'ees. D\'ecrivons
bri\`evement ces notions en suivant \cite{DerksenWeymanZelevinsky08}.
Soit en effet $Q$ un carquois fini. Notons $\hat{\C Q}$ {\em l'alg\`ebre
des chemins compl\'et\'ee}, c'est-\`a-dire la compl\'etion
de $\C Q$ par rapport \`a l'id\'eal bilat\`ere engendr\'e
par les fl\`eches de $Q$. L'espace $\C Q$ admet donc
une base topologique form\'ee de tous les chemins de~$Q$.
{\em L'homologie de Hochschild continue $\HH_0(\widehat{\C Q})$}
est le compl\'et\'e de l'espace quotient de $\C Q$ par
le sous-espace $[\C Q, \C Q]$ engendr\'e par tous les
commutateurs. Il admet une base topologique form\'ee
de tous les {\em cycles} de $Q$, c'est-\`a-dire les
orbites sous l'action du groupe cyclique $\Z/t\Z$
de chemins cycliques de longueur $t\geq 0$. Pour
chaque fl\`eche $\alpha$ de~$Q$, la {\em d\'eriv\'ee
cyclique} \cite{RotaSaganStein80}
est l'unique application lin\'eaire continue
\[
\partial_\alpha: \HH_0(\widehat{\C Q}) \to \widehat{\C Q}
\]
qui envoie la classe d'un chemin $p$ sur la somme $\sum vu$
prise sur toutes les d\'ecompositions $p=u\alpha v$
en des chemins $u$ et $v$ de longueur sup\'erieure ou \'egale
\`a z\'ero. Soit $W$ un {\em potentiel sur $Q$}, c'est-\`a-dire
un \'el\'ement de $\HH_0(\widehat{\C Q})$. {\em L'alg\`ebre
de Jacobi $\cp(Q,W)$} est la compl\'etion du quotient
de $\widehat{\C Q}$ par l'id\'eal bilat\`ere engendr\'e par
les d\'eriv\'ees cycliques $\partial_\alpha W$, o\`u $\alpha$
parcourt les fl\`eches de $Q$. Une {\em repr\'esentation d\'ecor\'ee
$(M,V)$} de $(Q,W)$ est form\'ee d'un module $M$ sur
$\cp(Q,W)$ et d'une famille d'espaces vectoriels $V_i$,
$i\in Q_0$, o\`u $M$ et les $V_i$ sont suppos\'es de dimension
finie sur $\C$. Par exemple, si $Q$ n'a pas de cycles
orient\'es (et donc $W=0$ et $\cp(Q,W)=\C Q$), toute
repr\'esentation d\'ecor\'ee $(M,V)$ fournit un
objet
\[
\pi(M)\oplus \bigoplus_{i\in Q_0} \Sigma \pi(V_i\ten_\C P_i)
\]
de la cat\'egorie amass\'ee (voir la section~\ref{ss:categorification-additive-categorie-amassee}).

Dans \cite{DerksenWeymanZelevinsky08} et \cite{DerksenWeymanZelevinsky09},
Derksen-Weyman-Zelevinsky  construisent et \'etudient l'op\'eration de
mutation pour les carquois \`a potentiel et leurs repr\'esentations d\'ecor\'ees.
Des difficult\'es techniques nombreuses et subtiles sont dues au
fait que cette op\'eration n'est ni fonctorielle ni d\'efinie partout.
Ceci est aussi la raison pour laquelle les repr\'esentations d\'ecor\'ees
ne forment pas, en g\'en\'eral, une cat\'egorie. N\'eanmoins,
la th\'eorie d\'evelopp\'ee par Derksen-Weyman-Zelevinsky est
assez proche de la cat\'egorification additive. Elle en diff\`ere
par le fait que les objets combinatoires centraux ne sont
plus les variables d'amas mais les $F$-polyn\^omes et $g$-vecteurs
introduits dans \cite{FominZelevinsky07} et qui sont peut-\^etre
encore plus fondamentaux que les variables d'amas.
Dans \cite{DerksenWeymanZelevinsky09}, Derksen-Weyman-Zelevinsky
appliquent leur th\'eorie en d\'emontrant de nombreuses conjectures
formul\'ees dans \cite{FominZelevinsky07}. Ils y parviennent
sous la seule hypoth\`ese que les alg\`ebres amass\'ees consid\'er\'ees
proviennent de bons carquois (non glac\'es), c'est-\`a-dire de matrices
antisym\'etriques quelconques, ce qui repr\'esente un
progr\`es remarquable par rapport aux approches pr\'ec\'edentes.
La construction de bases et la conjecture de positivit\'e restent
n\'eanmoins des probl\`emes compl\`etement ouverts dans cette
g\'en\'eralit\'e. Les id\'ees de Derksen-Weyman-Zelevinsky ont \'et\'e
li\'ees \`a la cat\'egorification additive au sens des
sections~\ref{ss:categorification-additive-preprojective}
et \ref{ss:categorification-additive-categorie-amassee} dans \cite{BuanIyamaReitenSmith08}
\cite{KellerYang09}~\cite{Amiot08a}. Un lien important avec
les surfaces \`a bord marqu\'ees est \'etabli dans \cite{Labardini09a}.

\bigskip\bigskip

 \subsection*{Remerciements}
 Je remercie 
 Caroline Gruson,
 Bernard Leclerc et
 Rached Mneimn\'e
 pour leurs conseils avis\'es sur une version
 ant\'erieure de ce texte.

\bigskip\bigskip\bigskip

\def\cprime{$'$} \def\cprime{$'$}
\providecommand{\bysame}{\leavevmode\hbox to3em{\hrulefill}\thinspace}

\end{document}